\documentclass[11pt,reqno]{amsproc}

\title[]{Remarks on a paper by Gavrilov: Grad-Shafranov equations, steady solutions of the three dimensional incompressible Euler equations with compactly supported velocities, and applications}

\author{Peter Constantin}
\address{Department of Mathematics, Princeton University, Princeton, NJ 08544}
\email{const@math.princeton.edu}
\author{Joonhyun La}
\address{Department of Mathematics, Princeton University, Princeton, NJ 08544}
\email{joonhyun@math.princeton.edu}
\author{Vlad Vicol}
\address{Courant Institute of Mathematical Sciences, New York University,
New York, NY, 10012}
\email{vicol@cims.nyu.edu}
\usepackage[margin=1in]{geometry}
\usepackage{amsmath, amsthm, amssymb}
\usepackage{times}
\usepackage{color}
\usepackage{hyperref}

\newcommand{\pa}{\partial}
\newcommand{\la}{\label}
\newcommand{\fr}{\frac}
\newcommand{\na}{\nabla}
\newcommand{\be}{\begin{equation}}
\newcommand{\ee}{\end{equation}}
\newcommand{\ba}{\begin{array}{l}}
\newcommand{\ea}{\end{array}}
\newcommand{\Rr}{{\mathbb R}}

\newcommand{\beg}{\begin}

\renewcommand{\div}{{\mbox{div}\,}}
\newcommand{\D}{\Delta}

\date{today}
\begin{document}
\begin{abstract}We describe a method to construct smooth and compactly supported solutions of 3D incompressible Euler equations and related models. The method is based on localizable Grad-Shafranov equations and is inspired by the recent result \cite{gav}. 
\end{abstract}

\noindent

\maketitle

\section{Introduction}
The three dimensional incompressible Euler equations are the basic equations of  mathematical fluid mechanics. The equations,
\be
\pa_t u + u\cdot\na u + \na p = 0,
\la{eulereq}
\ee
together with the incompressibility condition,
\be
\na\cdot u = 0,
\la{inc}
\ee
are four equations for the four unknown functions, velocities $u(x,t)\in \Rr^3$, and pressure $p(x,t)\in \Rr$, which depend on four independent variables, $x\in \Rr^3$, $t\in\Rr$.  
The pressure enforces the incompressibility condition, and thus obeys
\be
-\D p  = \na\cdot \left( u\cdot\na u\right ).
\la{preseq}
\ee
The Euler equations are conservative: smooth solutions preserve the energy
$\int_{\Rr^3}|u(x,t)|^2 dx$. If the pressure $p(x,t)$ belongs to $L^1(\Rr^3)$ at any instant of time, then
\be
\int_{\Rr^3}u_i(x,t)u_j(x,t)dx = -\delta_{ij}\int_{\Rr^3}p(x,t)dx.
\la{uijtriv}
\ee
This implies that the components of velocity are orthogonal and have equal norms in $L^2(\Rr^3)$, and the integral of pressure is negative. 
Moreover, if $u$ is Beltrami, i.e. if the curl of the velocity
\be
\omega = \na\times u
\la{omega}
\ee
is parallel to the velocity, and if $u\in L^2(\Rr^3)$, then $u$ must be identically zero (\cite{chaec}, \cite{nad}). In fact, Liouville theorems which assert the  vanishing of solutions which have constant behavior at infinity are often true for systems of the sort we are discussing.  In contrast, vortex rings are examples of solutions of the 3D Euler equations with compactly supported vorticity~\cite{FraenkelBurger74}. However, they have nonzero constant velocities at infinity.
Because of the Biot-Savart law
\be
u(x,t) = -\fr{1}{4\pi}\int_{\Rr^3}\fr{x-y}{|x-y|^3}\times \omega(y,t)dy,
\la{bs}
\ee 
if $\omega$ is compactly supported, it is hard to imagine that $u$ can also be compactly supported. In view of these considerations, the following result of Gavrilov (\cite{gav}) is surprising. 
\beg{thm}\la{gavthm}(Gavrilov, \cite{gav}) There exist nontrivial time independent solutions $u\in \left(C_0^{\infty}(\Rr^3)\right)^3$ of the three dimensional incompressible Euler equations.
\end{thm}
The purpose of this paper is to describe a proof of the result above, inspired by the original proof of Gavrilov, but starting from Grad-Shafranov equations, classical equations arising in the study of plasmas ~\cite{GradRubin58,Shafranov58} augmented by a {\em{localizability}} condition (see (\ref{constr})). This point of view yields a general method which can be applied to many other hydrodynamic equations, revealing a number of common features which seem rather universal. The 3D incompressible Euler equations result which extends Theorem~\ref{gavthm} is stated in Theorem~\ref{gsthm} below. An application providing  multiscale steady solutions which are locally smooth but belong only to H\"{o}lder classes $C^{\alpha}(\Omega)$ is given in Theorem \ref{turb}. For instance, locally smooth, compactly supported solutions can be constructed so that they belong to $L^2(\Omega)\cap C^{\fr{1}{3}}(\Omega)$ but not to any $C^{\alpha}(\Omega)$ with $\alpha>\fr{1}{3}$, they have vanishing local dissipation $u\cdot\na (\fr{|u|^2}{2} + p) = 0$, but have arbitrary large $\||\na u| |u|^2\|_{L^{\infty}(\Omega)}$. These solutions conserve energy, as they are stationary in time, and they have the regularity of the dissipative solutions recently constructed in connection with the Onsager conjecture (see review papers~\cite{BV,DLSZ}). These examples conclude the analysis of Sections \ref{axi} and \ref{construct}.  We construct by the same method time independent solutions of the incompressible 2D Boussinesq system with compactly supported velocities and
temperature in Section \ref{bousec} (~Theorem~\ref{bousthm}). It is known that smooth compactly supported solutions of the incompressible porous medium equations are identically zero (\cite{Tarek}). In Section \ref{ipmsec} we construct by the present method stationary solutions of the incompressible porous medium equation with velocities and temperature supported in arbitrary non-vertical strips (~Theorem~\ref{ipmthm}). 
\section{Steady Axisymmetric Euler Equations}\la{axi}
The stationary 3D axisymmetric Euler equations are solved via the Grad-Shafranov ansatz
\be
u = \fr{1}{r}(\pa_z \psi) e_r 
-\fr{1}{r}(\pa_r \psi )e_z  + \fr{1}{r}F(\psi)e_{\phi}
\la{u}
\ee
where $\psi = \psi(r,z)$ is a smooth function of $r>0$, $z\in\Rr$, and the swirl $F$ is a smooth function of $\psi$ alone. It is known that  smooth compactly supported velocities solving stationary axisymmetric 3D Euler equations must vanish identically if the swirl $F$ vanishes ~\cite{JiuXin09}. Above  $e_r$, $e_z$, $e_\phi$ are the orthonormal basis of cylindrical coordinates $r,z, \phi$ with the orientation convention $e_r\times e_\phi = e_z$, $e_r\times e_z = -e_{\phi}$, $e_\phi\times e_z = e_r$.  
Note that $u$ is automatically divergence-free, 
\be
\div u = 0,
\la{divz}
\ee
and also that, by construction,
\be
u\cdot\na\psi = 0.
\la{unapsi}
\ee
The vorticity $\omega = \na\times u$ is given by
\be
\omega = -\fr{1}{r}(\pa_z\psi) F'(\psi) e_r + \fr{1}{r}(\pa_r\psi) F'(\psi) e_z +
\fr{\D^*\psi}{r} e_{\phi}
\la{om}
\ee
where $F' =\fr{dF}{d\psi}$ and the Grad-Shafranov operator $\D^*$ is
\be
\D^*\psi = \pa_r^2\psi -\fr{1}{r}\pa_r\psi + \pa_z^2\psi.
\la{gsdelta}
\ee
In view of (\ref{u}) and (\ref{om}), the vorticity  can be written as
\be
\omega = -F'(\psi) u + \fr{1}{r}\left(\D^*\psi + \fr{1}{2}(F^2)'\right)e_{\phi}.
\la{omeg}
\ee
As it is very well known, the steady Euler equations
\be
u\cdot\na u + \na p = 0
\la{steadyeu}
\ee
can be written as
\be
\omega\times u +\na \left(\fr{|u|^2}{2} +p\right) = 0,
\la{omegasteadyeul}
\ee
and therefore the axisymmetric Euler equations are solved if $\psi$ solves the Grad-Shafranov equation (\cite{GradRubin58,Shafranov58}) 
\be
-\D^* \psi = \pa_{\psi}\left(\fr{F^2}{2} + r^2P\right )
\la{gs}
\ee
where the function $P = P(\psi)$ represents the plasma pressure:
\be
\omega \times u = \na P.
\la{omegaeq}
\ee
The analogy with the steady MHD equations  $u\leftrightarrow B$, $\omega\leftrightarrow J$ motivates the name. Both the swirl $F$ and the plasma pressure $P$ are arbitrary functions of $\psi$.  The plasma pressure and the hydrodynamic pressure are related via
\be
p + \fr{|u|^2}{2} + P = \mbox{constant}.
\la{pP}
\ee
The constant should be time independent if we are studying time independent solutions, and it may be taken without loss of generality to be zero.

In a remarkable paper, \cite{gav}, Gavrilov showed that smooth compactly supported velocities (\ref{u}) are possible. The construction of \cite{gav} is explicit but somewhat obscure, as it starts not from the Grad-Shafranov ansatz, but from a similar ansatz in terms of the hydrodynamic pressure $p$. The solution obtained by Gavrilov ends up though by having the hydrodynamic pressure proportional to the stream function $\psi$. Gavrilov's simple but important insight in 
\cite{gav} is that, if 
\be
u\cdot\na p = 0,
\la{crup}
\ee 
then, together with a solution $u,p$ of (\ref{steadyeu}, \ref{divz}), any function
\be
\widetilde{u} = \phi(p)u
\la{tildu}
\ee
with $\phi$ smooth is again a solution of (\ref{steadyeu}, \ref{divz}) with pressure given by
\be
\na{\widetilde p} = \phi^2(p)\na p.
\la{tildep}
\ee 
This can be used to localize solutions. In his construction Gavrilov obtained solutions $u$ defined in the neighborhood of a circle, obeying the Euler equations  near the circle, and having a relationship $|u|^2 = 3p$ between the velocity magnitude and the hydrodynamic pressure.

This motivates us to consider the overdetermined system formed by the Grad-Shafranov equation for $\psi$ (\ref{gs}) coupled with the requirement
\be
\fr{|u|^2}{2} = A(\psi).
\la{constr}
\ee
This requirement is the constraint of {\em{localizability}} of the Grad-Shafranov equation, and it severely curtails the freedom of choice of functions $F$ and $P$. This localizability is in fact the essence and the novelty of the 
method.  Without this constraint many solutions (\ref{u}) with $\psi$ solving the Grad-Shafranov equation (\ref{gs}) exist, including explicit ones (\cite{Sovolev}), but they cannot be localized in space.

The method we are describing consists thus in seeking functions $F,P, A$ of $\psi$ such that the system 
\be
\left\{
\ba
-\D^*\psi = \pa_{\psi}\left(\fr{1}{2}F^2(\psi) + r^2P(\psi)\right),\\
|\na\psi|^2 + F^2(\psi) = 2r^2A(\psi),
\ea
\la{sys}
\right.
\ee
is solved. Then the function $u$ given in the ansatz (\ref{u}), and the pressure
\be
p = -P(\psi) - A(\psi)
\la{p}
\ee
together satisfy the steady 3D Euler equations (\ref{steadyeu}, \ref{divz}), and are {\em{localizable}}, meaning that (\ref{constr}) is satisfied. It is important to observe that it is enough to find smooth functions $F, P, A$ of $\psi$ and a smooth function $\psi$ in an open set. This open set need not be simply connected, but once $u$ and $p$ are found using this construction, any function $\phi(p)u$ is again a solution of steady Euler equations, and it is sometimes possible to extend this solution to the whole space. 

\beg{rem}
Note that the functions $\rho = \sqrt{r^2+z^2}$, $\rho =  \fr{r^2}{2}$ and $\rho =z$ are all in the kernel of $\D^*$, i.e.
\be
\D^*\rho = 0,
\la{rhoeqs}
 \ee
and therefore, for $\rho = \sqrt{r^2+z^2}$ and $\rho =z$ the function $\psi = f(\rho)$ solves the vacuum ($P=0$) Grad-Shafranov equation if $f$ solves the ODE
\be
-f''(\rho) = FF'(f(\rho)) 
\la{feq}
\ee
because $|\na\rho |= 1$ in both cases. A function $\psi = f\left(\fr{r^2}{2}\right)$ solves the Grad-Shafranov equation with $F= {\mbox{constant}}$ if $f$ solves the ODE
\be
- f''(\rho) = P'(f(\rho)).
\la{feqr}
\ee
All these ODEs can be integrated (multiplying by $f'(\rho)$) but the solutions $\psi$ cannot be compactly supported in $\Rr^3$.

The two dimensional Euler equations have steady solutions with compactly supported velocities. Indeed, any smooth radial stream function produces a steady solution of the Euler equations in 2D, and if it is compactly supported in $\Rr^2\setminus \{0\}$ then the associated velocity is smooth and compactly supported.

\end{rem}

\section{Construction}\la{construct}
The construction of solutions of (\ref{sys}) starts with a hodograph transformation. We seek functions $U(r,\psi)$ and $V(r,\psi)$ defined in an open set in the $(r,\psi)$ plane and a smooth function $\psi(r, z)$ defined in an open set of the $(r,z)$ plane such that the equations
\be
\pa_r\psi(r,z) = U(r, \psi(r,z)),
\la{psireq}
\ee
\be
\pa_z\psi(r,z) = V(r, \psi(r,z))
\la{psizeq}
\ee
are satisfied. This clearly requires the compatibility
\be
V\pa_\psi U = U\pa_\psi V + \pa_r V.
\la{compa}
\ee
Once the compatibility is satisfied then the solution $\psi$ exists locally (in simply connected components).
The system (\ref{sys}) becomes
\be
\left\{
\ba
\pa_rU + U\pa_{\psi}U + V\pa_{\psi}V - \fr{1}{r}U = -F\pa_\psi F -r^2 \pa_\psi P\\
U^2+ V^2 +F^2  = 2r^2A.
\ea
\right.
\la{sysuv}
\ee
We traded a system of two equations in two independent variables ($r,z$) of total degree three, (\ref{sys}), for a system of three first order equations (\ref{compa}, \ref{sysuv}) in two independent variables ($r, \psi$). We integrate this locally. We start by noticing that the first equation of (\ref{sysuv}) is
\be
\pa_rU -\fr{1}{r}U + \fr{1}{2}\pa_{\psi}\left( U^2+ V^2 +F^2\right) = -r^2\pa_{\psi} P,
\la{newfirst}
\ee
which, in view of the second equation in (\ref{sysuv}), becomes
\be
\pa_rU -\fr{U}{r} = -r^2\pa_{\psi}(A+P),
\la{UAP}
\ee
and, using (\ref{p}) we see that
\be
\pa_{\psi} p = \fr{1}{r}\pa_r\left(\fr{U}{r}\right), 
\la{Up}
\ee
which then can be used to determine $p$ from knowledge of $U$. 
We observe that in order to have $p=p(\psi)$ a function of $\psi$ alone, from (\ref{Up}) we have to have
\be
U = r^3M(\psi) + rN(\psi).
\la{upsi}
\ee
for some functions $M$, $N$ of $\psi$.  Let us denote 
\be
Q_2(r,\psi) = 2r^2A(\psi)-F^2(\psi), 
\la{q2}
\ee
\be
Q_3(r,\psi) = r^3M(\psi) + rN(\psi),
\la{q3}
\ee
and
\be
Q_6(r,\psi) = Q_2(r,\psi) - (Q_3(r,\psi))^2
\la{q6}
\ee
polynomials of degree 2, 3 and 6 in $r$ with smooth and yet unknown coefficients depending only on $\psi$. We note that, in view of (\ref{upsi}), 
\be
U = Q_3,
\la{uq3}
\ee
and that the second equation in (\ref{sysuv}) yields 
\be 
V^2 = Q_6.
\la{vsquareq}
\ee
Multiplying (\ref{compa}) by $V$ results in 
\be
\pa_r Q_6 + Q_3\pa_\psi Q_6  - 2(\pa_\psi Q_3) Q_6 = 0.
\la{final}
\ee
Identifying coefficients in the 9th order polynomial equation (\ref{final}) we observe that only odd powers appear, the equations for powers 9 and 7 are tautological, and the remaining three equations become the ODE system
\be
\left\{
\ba
2(N^2-2A) + N(F^2)' -2N'F^2 =0,\\
8MN + M(F^2)'-2M'F^2 + N(N^2-2A)'-2N'(N^2-2A) = 0,
\ea
\right.
\la{odesy}
\ee
coupled with a separate equation for $A, M$, not involving $N,F$,
\be
A'M -3M^2 -2AM' = 0.
\la{am}
\ee
We recall that $' =\fr{d}{d\psi}$. This system of 3 first order ODEs with four unknown functions is equivalent to the compatibility relation (\ref{compa}). There is of course room to design solutions. 
Let us denote
\be
\alpha = \fr{F^2}{2A}
\la{alpha}
\ee
and
\be
\beta = - \fr{N}{M},
\la{beta}
\ee
assuming that $A,M$ have been chosen satisfying (\ref{am}). Then (\ref{odesy}) can be written as
\be
\alpha' = \fr{2}{M}\fr{1}{\alpha-\beta} + \fr{M}{A}(\beta- 3\alpha),
\la{alphaeq}
\ee
and
\be
\beta' = \fr{1}{M}\fr{1}{\alpha-\beta}.
\la{betaeq}
\ee
In order to localize the sought solution $u$ in $(r,z)$ space we need the pressure $p$ to take a value at a point $(r_0,z_0)$ which is strictly separated from all the values it takes on a circle in $(r,z)$ around that point. We seek then conditions which result in a strict local minimum for the function $\psi$ at the chosen point $(r_0,z_0)$, and then a similar behavior for the resulting $p$. Without loss of generality we may take this local minimum value of $\psi$ to be zero. We are lead thus to seek solutions  of the ODEs (\ref{am}, \ref{alphaeq}, \ref{betaeq})) 
in a small neighborhood of zero. Because $U$ and $V$ represent derivatives of $\psi$ we are lead to the requirement that the polynomials $Q_3$ and $Q_6$ both vanish at the point $(r_0, 0)$ in the  $(r,\psi)$ plane. This implies both $Q_3(r_0,0)=0$ and $Q_2(r_0,0)=0$. Note that
\be
Q_3(r,\psi) = rM(\psi)\left (r^2- \beta(\psi)\right) 
\la{q3beta}
\ee
and
\be
Q_2(r,\psi) = 2A(\psi)\left (r^2 -\alpha(\psi)\right).
\la{q2alpha}
\ee
The ODEs (\ref{alphaeq},\ref{betaeq}) are singular at $\alpha = \beta$. In order to still have vanishing of $Q_2, Q_3$, we choose solutions with $M(0)\neq 0$ and $A(0)=0$. The general solution of (\ref{am}) with given nonvanishing $M(t)$ and $A(0)= 0$ is
\be
A(t) = 3M^2(t)\int_0^t\fr{1}{M(s)}ds.
\la{atmt}
\ee
We choose the simplest one, the particular solution of (\ref{am})
\be
A(t) = 3mt, \quad M(t)=m
\la{amt}
\ee
with $m>0$ a constant. We wrote $t$ instead of $\psi$ for the one dimensional independent variable in order to avoid confusion with the sought function of $(r,z)$.

The initial datum for $\beta$ is determined by the choice of the location $r_0$. We take thus $\beta_0 = \sqrt{r_0}$ and 
\be
\beta(0) = \beta_0.
\la{betazero}
\ee
In view of the vanishing of $A$ at zero, in order for the equation (\ref{alphaeq}) to be nonsingular for small $t$ we need then to require  
\be
\alpha_0 = \fr{\beta_0}{3}.
\la{alphazero}
\ee
It is possible to locally solve the system (\ref{alphaeq}, \ref{betaeq}) with this choice of $A, M$. 
The solution of (\ref{alphaeq}, \ref{betaeq}) with coefficients (\ref{amt}) and initial data (\ref{betazero}, \ref{alphazero}) is most transparent if we write
\be
z(t) = \beta(t)- 3\alpha(t),
\la{zt}
\ee
with initial datum $z(0) =0$, and
\be
\zeta(t) = \fr{1}{\beta(t) - \alpha(t)},
\la{zeta}
\ee
with initial datum $\zeta(0) = \fr{1}{2\alpha_0} =\zeta_0$. The system (\ref{alphaeq}, \ref{betaeq}) becomes
\be
z' = -\fr{z}{t} + \fr{5}{m}\zeta
\la{z'}
\ee
coupled to
\be
\zeta' = \fr{z}{3t}\zeta^2 -\fr{1}{m}\zeta^3.
\la{zeta'}
\ee
We solve this with the ansatz
\be
z(t) = \sum_{j=1}^{\infty}z_jt^j, \quad \quad \zeta(t) = \sum_{j=0}^{\infty}\zeta_j t^j
\la{an}
\ee
and obtain from the first equation
\be
z_j = \fr{5}{m(j+1)}\zeta_{j-1},\quad j\ge 1.
\la{zj}
\ee
Equating coefficients of $t^{j-1}$ in the second equation we have
\be
j\zeta_j = \fr{1}{3}\left (z\zeta^2\right)_{j}  -\fr{1}{m}\left( \zeta^3\right)_{j-1}
\la{coeffj}
\ee
where $(f(t))_j$ means the coefficient of $t^j$ in the expansion of the function in power series. Using (\ref{zj}) it follows that
\be
\zeta_j = \fr{5}{3mj(j+1)}\zeta_0^2 \zeta_{j-1} + C_j(\zeta)
\la{zetaj}
\ee
where $C_j$ is a cubic convolution term with bounded coefficients and depending on the earlier variables $\zeta_k$, $k\le j-1$. 
This proves that the solution exists for short time $t$ and is analytic.
Consequently, $\alpha$ and $\beta$ exist and are analytic. In order to remove unnecessary parameters, let us rescale
\be
\phi = \fr{\psi}{m\beta_0^2},
\la{phi}
\ee
\be
r= \sqrt{\beta_0}(1 + x),
\la{x}
\ee
\be
z = y\sqrt{\beta_0},
\la{y}
\ee
\be
a(\phi)  = \fr{\alpha (\psi)}{\beta_0},
\la{anew}
\ee
\be
b(\phi) = \fr{\beta(\psi)}{\beta_0},
\la{bnew}
\ee
and define, after rescaling 
\be
P_3(x,\phi) = \fr{Q_3(r,\psi)}{m\beta_0^{\fr{3}{3}}} = (1+x)[(1+x)^2-b(\phi)]
\la{p3}
\ee
and
\be
P_2(x,\phi) =\fr{Q_2(r,\psi)}{m^2\beta_0^3} = 6\phi[(1+x)^2 -a(\phi)],
\la{p2}
\ee
resulting in
\be
P_6(x,\phi) = \fr{Q(r,\psi)}{m^2\beta_0^3} =  6\phi[(1+x)^2 -a(\phi)] - \left\{
(1+x)[(1+x)^2-b(\phi)]\right\}^2.
\la{p6}
\ee
These rescalings are natural for the Euler equations, they only look a bit unusual with our choice of constants. The length scale is $\ell = \sqrt{\beta_0}$ and the time scale is $\tau = (m\sqrt{\beta_0})^{-1}$. 

The ODEs solved by $a$ and $b$ are the rescaled (\ref{alphaeq}, \ref{betaeq}) with the rescaled (\ref{amt}):
\be
\fr{da}{d\phi} = \fr{2}{a-b} + \fr{1}{3\phi}(b-3a), \quad\quad  a(0) = \fr{1}{3},
\la{aode}
\ee
and
\be
\fr{db}{d\phi} = \fr{1}{a-b},\quad\quad b(0)=1.
\la{bode}
\ee
By the previous argument using (\ref{z'}, \ref{zeta'}) with $m=1, \alpha_0 =\fr{1}{3}$, their solutions are defined and analytic a small interval 
\be
\phi \in I = [-\epsilon, \epsilon].
\la{i}
\ee
The functions $a(\phi), b(\phi), Q_2(x,\phi), Q_3(x,\phi), Q_6(x,\phi)$ are well defined and analytic for any $x$, and $\phi \in I$, for small and fixed $\epsilon>0$. Moreover the compatibility equation (\ref{final}) is satisfied there
\be
\pa_x P_6  + P_3 \pa_{\phi} P_6 - 2 (\pa_\phi P_3)P_6 = 0.
\la{peq}
\ee

The equations  (\ref{psireq}), (\ref{psizeq}) become
\be
\pa_{x}\phi(x,y) = (1+x)[(1+x)^2- b(\phi(x,z))] = P_3(x,\phi(x,z))
\la{paxphi}
\ee
and
\be
(\pa_y\phi(x,y))^2 =  6\phi[(1+x)^2 -a(\phi(x,y))] - \left\{(1+x)[(1+x)^2-b(\phi(x,y))]\right\}^2 = P_6(x,\phi(x,y)).
\la{payphi}
\ee
These need to be solved in the neighborhood of $x=0$, $y=0$, and to yield a $C^1$ function  $\phi$ with $\phi(0,0) =0$, and $\phi> 0$ locally near $(0,0)$. There are no more parameters in the problem, and the equations are nondimensional.

We are thus interested in the open set $D$ in the $(x,\phi)$ plane
\be
D = \{(x,\phi)\left |\right.\; \phi \in I, P_6(x,\phi)> 0\}
\la{d}
\ee 
Note that $(0,0)$ is in the boundary of this set. Also, if $\phi(x,y)$ is differentiable in $x$ and solves (\ref{paxphi}) alone, irrespective of its initial datum, then
\be
\fr{d}{dx}P_6(x,\phi(x,y)) = 2(\pa_{\phi}P_3(x,\phi(x,y)))P_6(x,\phi(x,y))=
-2(1+x)b'(\phi(x,y))P_6(x,\phi(x,y)).
\la{psixchar}
\ee
Consequently, if some $(x_0, \phi(x_0, y_0))\in D$, then the locally $(x,\phi(x,y_0))$ must be in $D$.  Furthermore, if $P_6(x_0,\phi_0) = 0$, then on the solution it stays zero: the boundary $\pa D$ is characteristic. 

Let us consider the boundary of $D$ more closely. By the implicit function theorem, in view of the fact that
\be
\pa_\phi P_6(0,0) = 4\neq 0,
\la{pasixprimezero}
\ee
we have the existence of a smooth function $\delta: [-\epsilon, \epsilon]\to \Rr$, 
\be
x\mapsto \delta(x), 
\la{deltax}
\ee
with $\delta(0)=0$ and satisfying on $I\times I$
\be
P_6(x,\phi) = 0 \Leftrightarrow \phi = \delta(x).
\la{implicit}
\ee
We might need to shrink the size of $I$, but we use the same $\epsilon$ as in the definition (\ref{i}) of $I$, without loss of generality, in order not to clutter the notation.
From the definition of $P_6$ and the fact that $a(0) = \fr{1}{3}$ it follows that in fact 
\be
\delta(x) \ge 0, \quad\mbox{for}\; x\in I.
\la{deltapos}
\ee
Also, computing with the aid of (\ref{peq}) we see that
\be
\fr{d\delta}{dx}(x) = - \fr{\pa_xP_6 (x,\delta)}{P'_6(x,\delta)} = P_3(x,\delta) - 2 P_3'(x,\delta)\fr{P_6(x,\delta)}{P_6'(x,\delta)},
\la{deltaprime}
\ee 
(where $' = \fr{d}{d\phi}$) and, using (\ref{implicit}) we deduce, importantly
\be
\fr{d\delta}{dx}(x) = P_3(x,\delta(x)).
\la{deltaU}
\ee
Consequently
\be
\fr{d\delta}{dx}(0) = 0,
\la{vanishdeltaprime}
\ee
and
\be
\fr{d^2\delta}{dx^2}(0) = \pa _x P_3(0,0) = 2.
\la{convexity} 
\ee
Therefore  $\delta(0) =0 $ is the strict local minimum of $\delta(x)$, and locally $D$ is the supergraph
\be
D = \{(x,\phi)\in I\times I\left |\right. \phi>\delta(x)\},
\la{Dgraph}
\ee
which means that if $(x_0,\phi_0)\in D$ then $(x_0,\phi)\in D$ for $\phi\in I$, $\phi>\phi_0$.  We define now for $(x,\phi)\in I\times I \cap D$ and $y\neq 0$,

\be
V(x,y,\phi) =
 \left\{
\ba
\sqrt{P_6(x,\phi)}, \quad\quad\quad  {\mbox {when}}\; y>0\\
-\sqrt{P_6(x,\phi)}, \quad\quad {\mbox{when}}\; y<0,
\ea
\right.
\la{vdef}
\ee
i.e.,
\be
V(x,y, \phi) ={\mbox{sgn}}\,(y)\sqrt{P_6(x,\phi)},
\la{vdefsign}
\ee
and we recall that, independently of $y$, 
\be
U(x,\phi) = P_3(x,\phi).
\la{udef}
\ee
We set
\be
\phi(x,0) = \delta(x) 
\la{id}
\ee
for all $x\in I$. We remark that, in view of (\ref{deltaU}), we have
\be
\delta(x) = \int_0^xP_3(\xi,\delta(\xi))d\xi.
\la{deltaint}
\ee
Then, for $y>0$ and $y<0$, (each case separately), we define
\be
\phi(x,y) = \int_0^y V(0,\eta, \phi(0,\eta))d\eta + \int_0^x P_3(\xi, \phi(\xi, y))d\xi
\la{pathtwo}
\ee
where $\phi(0,\eta)$ is the positive smooth solution of
\be
\fr{d}{dy}\phi(0,y) = V(0,y,\phi(0,y))
\la{left}
\ee
with initial data $\phi(0,0)=0$.  This requires a discussion because this is a non-Lipschitz ODE with infinitely many solutions, including an identically zero one. 
The function 
\be
\pi(\phi) = P_6(0,\phi) =6\phi(1-a(\phi)) - (1-b(\phi))^2\la{piphi}
\ee
is positive for $\phi>0$ and vanishes like $4\phi$ at $\phi =0$. The primitive
\be
Y(\phi) = \int_0^{\phi}\fr{dt}{\sqrt{\pi(t)}}
\la{yphi}
\ee
is positive, increasing, vanishes at $0$, is smooth on $(0, \epsilon]$, and has
\be
Y'(0) = \infty.
\la{yprimezero}
\ee
Because of monotonicity, $Y^{-1}$ is locally defined, and we put 
\be
\phi(0,y) = Y^{-1}(y)
\la{phizeroy}
\ee
for $y\ge 0$.
Note that 
\be
\fr{d}{dy}\phi(0, y) = V(0,y, \phi(0,y))
\la{eqzeroy}
\ee
holds fror all $0\le y\le \epsilon$. A similar construction is done for $y\le 0$. In fact we are constructing via (\ref{pathtwo}) two functions, $\phi_{\pm}(x,y)$, one for $y\ge 0$ and one for $y\le 0$, glued by (\ref{id}, \ref{deltaint}) and satisfying by construction
\be
\phi_{-}(x,-y) = \phi_+(x,y).
\la{evenphi}
\ee

Obviously  $P_3(x,\phi(x,y)) = \pa_x\phi(x,y)$, that is (\ref{paxphi}), is true by construction from (\ref{pathtwo}) and from (\ref{deltaU}). Now we claim that $\phi(x,y)\in D$. Indeed,
\be
\fr{d}{dx}\left(\phi(x,y) - \delta(x)\right ) = P_3(\phi(x,y))-P_3(\delta(x)) =
I(\left(\phi(x,y) - \delta(x)\right )
\la{phiind}
\ee
where $I=I(x,y) = \int_0^1P'_3(\delta(x)+ t(\phi(x,y))-P_3(\delta(x)))dt$ is bounded. The initial data is positive $\phi(0,y)-0>0$. Therefore $\phi(x,y)-\delta(x)>0$ and $\phi\in D$. (The same conclusion follows also from (\ref{psixchar}) which needs only (\ref{paxphi})). 

The vertical derivative is obtained differentiating the construction (\ref{pathtwo}) and we see that it is given by
\be
\pa_y\phi(x,y) = V(0,y,\phi(0,y)) + \int_0^xP'_3(\xi,\phi(\xi,y))\pa_y\phi(\xi,y)d\xi
\la{vertder}
\ee
Its $x$ derivative at fixed $y$ obeys
\be
\fr{d}{dx}Z(x)= P_3'(x,\phi(x,y))Z(x)
\la{zeq}
\ee
and this is also the equation obeyed by $V(x,y,\phi(x,y))$ in virtue of the fact that $\phi\in D$ and the compatibility equation (\ref{peq}). This again uses only (\ref{paxphi}). The values at $x=0$ equal both $V(0,y,\phi(0,y))$ and because $P'_3(x,\phi(x,y))$ is bounded it follows that
\be
\fr{d\phi(x,y)}{dy} = V(x,y,\phi(x,y))
\la{vertode}
\ee
for $y>0$, respectively $y<0$. 
We have thus the nontrivial solution $\phi(x,y)\in D$ for $y\neq 0$, and, from  the compatibility (\ref{peq}) we have the equations (\ref{paxphi}) and (\ref{payphi}) verified. 
We note that, in view of (\ref{deltaU}) we have that
\be
\phi(x,y) = \int_0^x P_3(\xi, \phi(\xi,0))d\xi + \int_0^y V(x,\eta,\phi(x,\eta))d\eta.
\la{pathone}
\ee
which is natural because once $\phi\in C^1$ has been constructed and (\ref{paxphi}) and (\ref{payphi}) have been verified, then the one form
\be
P_3(x,\phi(x,y))dx + V(x,y,\phi(x,y))dy
\la{form}
\ee
is closed. The function $\phi$ is $C^1$ near the origin, and the functions $V(x,y,\phi)$ and $P_3(x,\phi(x,y)$ are continuous. Successive differentiations of (\ref{paxphi}) and (\ref{payphi}) imply that $\phi$ is actually smooth.

In order to understand this better let us compare with the (simpler to compute) case in which we need to solve ODEs
\be
\left (\fr{d\phi(x,y)}{dy}\right )^2 = \phi(x,y) -\delta(x), \quad \phi(x,0) = \delta(x)
\la{odex}
\ee
In this case $P_6(x,\phi) = \phi -\delta(x)$, $D = \{\phi>\delta(x)\}$.  There
is no uniqueness, the boundary $\pa D$ is characteristic, and, in addition to
the solution $\phi(x,y)= \delta(x)$ we have the analytic solution 
$\phi(x, y) = \fr{y^2}{4} + \delta(x)$ and infinitely many discontinuous solutions
\be
\phi(x,y) = \left\{
\ba
\delta(x), \quad\quad {\mbox{for}}\; |x|<\epsilon_1,\\
\fr{y^2}{4} + \delta(x), \quad {\mbox{for}}\; \epsilon_1\le |x|<\epsilon
\ea
\right.
\la{exabad}
\ee
for any $\epsilon_1<\epsilon$. The case of the true $P_6(x,\phi)$ is similar.


We have proved therefore
\beg{thm}\la{gsthm} Let $\ell >0, \; \tau > 0\in \Rr$ be given. There exists $\epsilon>0$ and a function $\psi\in C^{\infty}(B)$, where $B = \{(r,z)\left |\right.\; |r-\ell|^2 + |z|^2 <\epsilon^2\ell^2\}$ satisfying  $\psi(\ell, 0) =0$, $\psi>0$ in $B$ and such that (\ref{sys}) holds with $A$, $P$ and $F^2 $ real analytic functions of $\psi$. The Grad-Shafranov equation (\ref{gs}) is solved pointwise and has  classical solutions in $B$. The associated velocity $u$ given by the Grad-Shafranov ansatz (\ref{u}) is H\"{o}lder continuous in $B$. The Euler equation (\ref{steadyeu}, \ref{divz}) holds weakly in $B$. The pressure is given by $p = \fr{1}{\ell\tau}\psi$. The vorticity is bounded, $\omega \in L^{\infty}(B)$ and (\ref{omegasteadyeul}) holds a.e. in $B$. 
\end{thm}
We wrote the theorem in a nondimensional form. In order to return to dimensional variables we translate in $z$, so that $z_0 = 0$ dilate in $r$ so that $r_0 = \ell = \sqrt{\beta_0}$ and choose $2m = \fr{1}{\ell\tau}$. 
We note that $F(\psi)$ vanishes like $\sqrt{\psi}$. Therefore, while the ansatz
(\ref{u}) gives a bounded swirl and a H\"{o}lder continuous velocity, the vorticity is not smooth. In fact, in view of (\ref{omeg}) the vorticity equals
\be
\omega (r,z) = - F'(\psi) u(r,z) + {\mbox{smooth}}.
\la{omegas}
\ee
Thus,  $\omega\in L^{\infty}(B)$, because $u$ vanishes to first order at $(\ell, 0)$, but the $r$ derivative of the $z$ component of vorticity is infinite there.

Once $\psi$ has been constructed so that it has a local minimum at $(\ell, 0)$, then $p$ has also a local minimum there, because, by (\ref{Up}),
\be
p = 2m\psi = \fr{\psi}{\ell\tau}
\la{ppsi}
\ee
is monotonic in $\psi$.  

Theorem \ref{gavthm} holds because the cutoff can be chosen so that the point $(\ell, 0)$ is omitted. By choosing a suitable cutoff function $\phi_\epsilon(p)$, the solution $\widetilde{u} = \phi_\epsilon(p) u$ is supported in the region $A = \{(r,z)\left|\right.\; \fr{1}{2}\ell^2\epsilon^2 <|r-\ell|^2 + |z|^2  <\epsilon^2\ell^2\}$.

An immediate consequence of Theorem~\ref{gsthm} is the following.

\beg{thm}\la{turb} Let $0\le \alpha<1$. In any domain $\Omega\subset\Rr^3$ there exist steady solutions of Euler equations belonging to $C^{\alpha}(\Omega)$ and vanishing outside a compact included in $\Omega$, but such that they do not belong to $C^{\beta}(\Omega)$ for $\beta>\alpha$. The solutions are locally smooth, meaning that for every $x\in \Omega$ there exists a neighborhood of $x$ where the solution is $C^{\infty}$. For any $\Gamma>0$, there exist such locally smooth steady solutions $u$ which belong to $L^2(\Omega)\cap C^{\fr{1}{3}}(\Omega)$, are supported in a compact subset of $\Omega$ and have
\[
\sup_{x\in\Omega}|\na u(x)||u(x)|^2 \ge \Gamma,
\]
while the local dissipation vanishes, i.e. $u\cdot \na(\fr{|u|^2}{2} + p) =0$ in the sense of distributions. There exist steady solutions which are locally smooth and whose Lagrangian trajectories have arbitrary linking numbers.
\end{thm}
\noindent{\bf{Proof.}} Taking a smooth template $u_B$ constructed in Theorem \ref{gsthm} and compactly supported in the toroidal shell $B$  in $\Rr^3$ 
\be
B = \{(r,z)\left| \right. \fr{1}{2} <(|r-1|^2 + z^2<1\},
\la{bone}
\ee
we write
\be
u(x) = \sum_{n=1}^\infty U_n R_n u_B\left (\fr{R^T_n(x-x_n)}{\ell_n}\right) = \sum_{n=1}^{\infty}u_n(x)
\la{uuB}
\ee
where the positive numbers $U_n>0$, $\ell_n>0$, the rotations $R_n\in O(3)$ with their transposed $R_n^T$, and the vectors $x_n\in \Rr^3$ are arbitrary.
We may take for instance an helicoidal curve of infinite length inside a relatively compact open subset of $\Omega$ and decorate it with small disjoint toroidal shells of scales $\ell_n$, centered at points $x_n$ on the curve and whose symmetry axes are tangent to it (determing thus $R_n$). Then the supports of the terms in the sum are mutually disjoint and $u$ given above is a steady solution of the incompressible Euler equations, because $u_B$ is one, and thus, by scaling and rotating,
\be
u_n\cdot\na u_n + \na p_n = 0.
\la{uneuler}
\ee
The supports of the gradients of pressure $\na p_n$ of the individual terms in the sum are mutually disjoint.
The norm of $u$ in $W^{m,p}$ is proportional to
\be
\|u\|_{W^{\alpha, p}(\Omega)}\sim \left (\sum_{n=1}^{\infty}U_n^p \ell_n^{(3-p\alpha)}\right)^{\fr{1}{p}}
\la{normuub}
\ee
for any $\alpha \ge 0$, $1\le p \le\infty$. In particular, we can find square integrable time independent solutions $u$ such that $u\in C^{\fr{1}{3}}(\Omega)$ but $u\notin C^{\beta}(\Omega)$ for $\beta>\fr{1}{3}$. The $L^{\infty}(\Omega)$ norms of $\pa^\alpha u$ are of the order $\sup_{n\ge 1}U_n\ell_n^{-\alpha}$.
 
\section{2D Boussinesq}\la{bousec}
The time independent 2D Boussinesq system  (\cite{bous}, \cite[Section 2.4]{Const17})
is
\be
\left\{
\ba
u\cdot\na u + \na p = \theta e_2,\\
u\cdot\na \theta = 0,\\
\na\cdot u = 0,
\ea
\right.
\la{bous}
\ee
in $\Rr^2$. The function $u= (u_1,u_2)$ is the incompressible velocity, $\theta$ is the temperature and $e_2$ is the direction of gravity.
We write this as
\be
\omega u^{\perp} = \na P + \theta e_2
\la{omegabous}
\ee
where $\omega = \pa_1 u_2-\pa_2 u_1$ is the vorticity and $u^{\perp} = -u_2e_1+ u_1e_2$. In terms of the stream function $\psi$ we have
\be
u = \na^{\perp}\psi,
\la{upsibous}
\ee
\be
\omega = \D\psi.
\la{omegapsi}
\ee
The equation (\ref{omegabous}) is therefore
\be
-(\D\psi)\na \psi = \na P + \theta e_2
\la{gsinter}
\ee
and choosing
\be
\theta(x) = \Theta(\psi(x))
\la{thetas}
\ee
and 
\be
P(x) = -x_2\Theta(\psi(x)) - G(\psi(x))
\la{P}
\ee
with $G(\psi)$ an arbitrary function, we see that the steady Boussinesq equation is satisfied if $\psi$ satisfies the corresponding Grad-Shafranov-like equation
\be
\D\psi = x_2\Theta'(\psi) + G'(\psi)
\la{gsbous}
\ee
where, like before $' = \fr{d}{d\psi}$. Identifying the hydrodynamic pressure, and requiring it to be a function of $\psi$ alone we see that
\be
p(\psi)  = G + x_2\Theta -\fr{1}{2}|\na\psi|^2.
\la{pgs}
\ee
Let us seek solutions of (\ref{gsbous}, \ref{pgs}) by hodograph transformation
\be
\pa_2\psi(x_1,x_2) = U(x_2,\psi(x_1,x_2))
\la{pa2psi}
\ee
and
\be
\pa_1\psi(x_1,x_2) = V(x_2, \psi(x_1,x_2)).
\la{pa1psi}
\ee
The compatibility equation is
\be
VU' = \pa_2V + UV'
\la{compabous}
\ee
and the Grad-Shafranov-like equation (\ref{gsbous}) becomes
\be
UU'+ \pa_2U + VV' =x_2\Theta' + G'.
\la{gsinteruv}
\ee 
Using (\ref{pgs}) we deduce that
\be
\pa_2 U(x_2, \psi) = p'(\psi).
\la{pa2Up}
\ee
This requires $U$ to be a linear polynomial
\be
U = p'(\psi)(x_2-\beta(\psi)).
\la{ubeta}
\ee
Defining the polynomial
\be
Q_1(x_2,\psi) = 2(x_2\Theta(\psi) + H(\psi))
\la{q1}
\ee
where
\be
H(\psi) = G(\psi) - p(\psi)
\la{hpsi}
\ee
we see that (\ref{pgs}) becomes
\be
V^2 = Q_1-U^2.
\la{vqu}
\ee
To summarize, after the transformation, we arrived at the equations (\ref{compabous}), (\ref{ubeta}) and (\ref{vqu}). Now (\ref{compabous}), after multiplication by $V$ becomes a polynomial identity,
\be
2V^2U' = \pa_2V^2 + U(V^2)'
\la{compauv}
\ee
which, using (\ref{vqu}) simplifies to read 
\be
2Q_1U'-UQ_1' + 2U\pa_2U -\pa_2Q_1 = 0.
\la{compauq}
\ee
Identifying coefficients in the ensuing ODEs we deduce
\be
\Theta = k(p')^2,
\la{thetap'}
\ee
with $k$ a constant, which of course we need to be nonzero,
\be
\Theta = -2p'\beta'(\beta\Theta + H),
\la{thetabeta'}
\ee
and
\be
(p')^2 -2\Theta p'\beta' + 2Hp'' -p'H' = 0.
\la{p'square}
\ee
Now we would like to localize, so we would need $U$ and $V$ to vanish at $(x_2^0,0)$. From (\ref{ubeta}) this requires either $p'(0)=0$ or $x_2^0 = \beta(0)$. If the latter, then the vanishing of $V$ requires the vanishing of $Q_1$ and that implies $\beta(0)\Theta(0) + H(0)=0$. But, in view of (\ref{thetabeta'}) this, in turn implies that $\Theta(0) = 0$ or $p'\beta'$ to be infinite at $0$.
If $p'\beta'$ is finite, then by (\ref{thetap'}) we are lead to $p'(0)= 0$. So, if we wish to localize we take 
\be
p'(0) = 0.
\la{pprimez}
\ee
We denote
\be
q= p'
\la{q}
\ee
and define $\alpha$ via
\be
H = \alpha q^2.
\la{alphah}
\ee
The equations for $\alpha$ and $\beta$ become
\be
q\alpha' = 1 + \fr{k^2}{\alpha + k\beta}
\la{qalpha'}
\ee
and
\be
q\beta' = -\fr{k}{2(\alpha + k\beta)}.
\la{qbeta'}
\ee
Introducing
\be
\gamma = \alpha + k\beta
\la{gamma}
\ee
we deduce
\be
q\gamma' = 1 + \fr{k^2}{2\gamma}
\la{qgamma'}
\ee
We integrate the latter, denoting
\be
q\fr{d}{d\psi} = \fr{d}{d\tau}
\la{qtau}
\ee
and we obtain
\be
\gamma - \fr{k^2}{2}\log \left(\gamma + \fr{k^2}{2}\right ) = \tau + C.
\la{gammasolved}
\ee
We choose
\be
\gamma(0) = 1
\la{gammaz}
\ee
by adjusting $C$ in (\ref{gammasolved}) appropriately. We then have from (\ref{gammasolved}), which we solve implicitly, a unique smooth function 
 \be
\gamma = \gamma(\tau),
\la{gammatau}
\ee
that is defined for small $\tau$ and is bounded away from zero. The equations (\ref{qalpha'}, \ref{qbeta'}) are 
\be
\fr{d\alpha}{d\tau} = 1 + \fr{k^2}{\gamma(\tau)}
\la{alphatau}
\ee
and
\be
\fr{d\beta}{d\tau} = -\fr{k}{2\gamma(\tau)}.
\la{betatau}
\ee
These yield functions
\be
\alpha(\tau), \beta(\tau)
\la{alphabeta}
\ee
which are smooth, finite for small $\tau$, and remain close to values $\alpha(0)$, $\beta(0)$ that are not restricted, except by our choice of $\gamma(0) = 1$,
\be
\alpha(0) + k\beta(0) = 1.
\la{idgamma}
\ee 
We abused notation by using the same letters $\alpha, \beta$ for the functions
of $\tau$ solving the explicit ODEs in $\tau$ (\ref{alphatau}, \ref{betatau}).
In fact a more precise notation would have been to use different letters, $\widetilde{\alpha}, \widetilde{\beta}$ and write $\alpha(\psi) = \widetilde{\alpha}(\tau), \beta(\psi) =\widetilde{\beta}(\tau)$ where $\tau(\psi)$ is determined by the defining relation (\ref{qtau}), 
\be
\fr{d\tau}{d\psi} = q^{-1}.
\la{dtaudpsi}
\ee
If we set
\be
q = \psi^{s},
\la{qpsi}
\ee
with $0<s<1$, and seek $\psi>0$, then 
\be
\tau  = \fr{1}{1-s}\psi^{1-s}.
\la{taupsi}
\ee
Returning to the equation (\ref{vqu}) for $V$,  we see that
\be
V^2 = q^2\left(-x_2^2 + 2x_2(k+\beta) + 2\alpha -\beta^2 \right).
\la{valphabeta}
\ee
Choosing for instance $x_2^0 = \beta(0)$ gives that the term in paranthesis
equals 2 at $\psi=0$. Thus, we have that
\be
V = \psi^s\sqrt{-x_2^2 + 2x_2(k+\beta) + 2\alpha-\beta^2}
\la{vxalphabeta}
\ee
and 
\be
U = \psi^s(x_2-\beta)
\la{uxbeta}
\ee
are H\"{o}lder continuous. For instance, if $s=\fr{1}{2}$ we obtain equations for $\sqrt{\psi}$ with coefficients which depend smoothly on $\sqrt{\psi}$. These can be integrated locally easily with the result that
\be
\psi(x^0) = 0
\la{psizero}
\ee
and $\psi\in C^{\infty}(B)$, where $B$ is a small ball around $x^0$. 
Notice that
\be
p = \fr{1}{1+s}\psi^{1+s}
\la{ppsix}
\ee
yields then a regular function with local minimum at $x^0$ and therefore we can localize the system. 

By taking $\widetilde{u} = \phi(p) u$, $\widetilde{\theta} = \phi(p)\theta$ and $\na\widetilde p = \phi^2(p)\na p$ with appropriate smooth cutoff $\phi$ we have obtained thus the following result.
\beg{thm}\la{bousthm} There exist smooth nontrivial time independent solutions of (\ref{bous}) such that $(u,\theta)$ are compactly supported.
\end{thm}
\section{Incompressible Porous Medium Equation}\la{ipmsec}
The 2D time independent inviscid porous medium (IPM) equation (\cite{ipm}, \cite[Section 2.3]{Const17}) is
\be
\left\{
\ba
u\cdot \na \theta = 0,\\
\na\cdot u = 0,\\
u = \theta e_2 + \na p,
\ea
\right.
\la{ipm}
\ee
where $e_2$ is the direction of gravity. We write  $u= \na^{\perp}\psi$ and require 
\be
p = p(\psi)
\la{ppsipm}
\ee
and
\be
\theta = \Theta(\psi).
\la{thetaipm}
\ee
The equations are
\be
\left\{
\ba
-\pa_y \psi = p'\pa_x\psi,\\
\pa_x \psi = p'\pa_y\psi + \Theta
\ea
\right.
\la{ipmpsi}
\ee
Consequently
\be
\pa_x\psi = \fr{\Theta(\psi)}{1 + (p'(\psi))^2} = U(\psi)
\la{paxpsipm}
\ee
and
\be 
\pa_y \psi = -\fr{\Theta(\psi)p'(\psi)}{1 + (p'(\psi))^2} = V(\psi).
\la{paypsipm}
\ee
In view of the fact that $\pa_x\psi  = U(\psi)$ and $\pa_y\psi = -p'(\psi)U(\psi)$ the compatibility equation is 
\be
p'(\psi) = k
\la{compaipm}
\ee
where $k$ is a constant. 
Because
\be
\pa_y\psi = -k\pa_x\psi
\la{linbal}
\ee
then the solution must depend on a single variable
\be
\psi (x,y) = f(x-ky).
\la{psif}
\ee
The equations (\ref{paxpsipm}, \ref{paypsipm}) reduce to
\be
f'(z) = \fr{1}{1+k^2}\Theta(f(z))
\la{f'}
\ee 
where 
\be
z = x-ky\la{z}.
\ee
Because we seek to localize, we look for functions $\Theta(\psi)$ which vanish at $\psi=0$. Taking for instance
\be
\Theta(\psi) = \psi^s
\la{thetapsipm}
\ee
with $0<s<1$ we obtain
\be
f(z) = \left(\fr{(1-s)(z-z_0)}{1+k^2}\right)^{\fr{1}{1-s}}
\la{fz}
\ee
and then we have 
\be
\psi(x,y) = \left(\fr{(1-s)(x-x_0-k(y-y_0))}{1+k^2}\right)^{\fr{1}{1-s}},
\la{pxixyipm}
\ee
and
\be
p(x,y)= k \left(\fr{(1-s)(x-x_0- k(y-y_0))}{1+k^2}\right)^{\fr{1}{1-s}}.
\la{pxyipm}
\ee
If $s=\fr{1}{2}$ we obtain parabolas, $\psi = c_kz^2, p = k\psi$  which can be localized in the $z$ variable.

By setting $\widetilde{u} = \phi(p) u$, $\widetilde{\theta} = \phi (p) \theta$ and $\na\widetilde{p} = \phi(p)\na p$  with appropriate smooth cutoff $\phi$, we have thus proved the following result.

\beg{thm}\la{ipmthm} Let $S$ be any strip in $\Rr^2$ of finite width, and whose
direction is not parallel to the direction of gravity $e_2$. There exist smooth time independent solutions of the incompressible porous medium equation (\ref{ipm}) for which the velocities $u$ and the temperature $\theta$ vanish outside the strip $S$.
\end{thm}

\noindent{\bf{Acknowledgments.}} The work of PC was partially supported by NSF grant DMS-1713985 and by the Simons Center for Hidden Symmetries and Fusion Energy.  JL was partially supported by a Samsung Fellowship. VV was partially supported by NSF grant DMS-1652134.

\end{document}